\documentclass[11pt]{article}   	                 
\usepackage{amsmath,amsthm, amssymb,amsfonts,amscd,mathabx}
\usepackage[pdftex]{graphicx}                 
\usepackage{wrapfig}                               
\usepackage[font=footnotesize]{caption}                     
\usepackage{float}

\usepackage[colorlinks=true,linkcolor=blue,citecolor=red]{hyperref}
\numberwithin{equation}{section}
\newtheorem{theorem}{Theorem}[section]

\newtheorem{proposition}[theorem]{Proposition}
\newtheorem{definition}[theorem]{Definition}

\newtheorem{remark}[theorem]{Remark}
\newtheorem{result}[theorem]{Result}
{\begin{trivlist}\item[]\textbf{Proof#1 }}%
{\qed\end{trivlist}}

\newenvironment{Acknowledgment}%
 {\begin{trivlist}\item[]\textbf{Acknowledgments }}{\end{trivlist}}

\newcommand{\R}{\mathbb{R}}
\newcommand{\C}{\mathbb{C}}
\newcommand{\N}{\mathbb{N}}

\newcommand{\Rmnum}[1]{\uppercase\expandafter{\romannumeral #1\relax}}

\newcommand{\rmo}{\mathrm{o}}

\newcommand{\rme}{\mathrm{e}}
\newcommand{\rmi}{\mathrm{i}}

\renewcommand{\Re}{\mathrm{Re}\,}
\renewcommand{\Im}{\mathrm{Im}\,}

\renewcommand{\leq}{\leqslant}
\renewcommand{\geq}{\geqslant}

\def\eps{\varepsilon}

\def\v{\mathbf{v}}

\newfam\bifam
\font\tenbi=cmmib10 scaled \magstep1 \font\sevenbi=cmmib10 at 11pt
\font\fivebi=cmmib10 at 6pt \textfont\bifam = \tenbi
\scriptfont\bifam = \sevenbi \scriptscriptfont\bifam= \fivebi

\begin{document}

\begin{center}
{\fontsize{14}{14}\fontfamily{cmr}\fontseries{b}\selectfont{
Pulled fronts are not (just) pulled}}\\[0.2in]
Montie Avery$^{1}$, Matt Holzer$^{2}$, and  Arnd Scheel$^{3}$\\[0.1in]
\textit{\footnotesize
$^1$Department of Mathematics, Emory University\\[0.05in]
$^2$ Department of Mathematical Sciences, George Mason University, 4400 University Drive,  Fairfax, VA, 22030, USA \\[0.05in]
$^3$School of Mathematics, University of Minnesota, Minneapolis, 206 Church St SE, MN 55414, USA
}
\end{center}

\begin{abstract}
\noindent
Front propagation into unstable states is often determined by the linearization, that is, propagation speeds agree with predictions from the linearized equation at the unstable state. The leading edge behavior is then a Gaussian tail propagating with the linear spreading speed. Fronts following this leading edge are commonly referred to as pulled fronts, alluding to the idea that they are ``pulled'' by this leading-edge Gaussian tail. We describe here a class of examples that exhibits how these leading-order effects do not completely describe the dynamics in the wake of the front. In fact, leading edge behavior predicts at most two possible invasion scenarios, associated with positive and negative amplitudes of the Gaussian tail, but our examples exhibit three or more invasion fronts with different states in the wake. The resulting invasion process therefore leaves behind a state that is not solely  determined by the leading edge, and thus not just pulled by the Gaussian tail.
\end{abstract}

\begin{Acknowledgment}

 The authors gratefully acknowledge partial support through grants  NSF DMS-2510541 and DMS-2202714 (M.A), NSF DMS-2406623 (M.H.) , and
 NSF DMS-2205663 and DMS-2506837 (A.S.).
\end{Acknowledgment}

%
%


\section{Introduction}

Localized perturbations of unstable states in spatially extended systems both grow in time and spread spatially.
Generally speaking, instabilities in large systems lead to very high-dimensional unstable manifolds and the evolution of the instability is difficult to predict. Initial perturbations that are localized in space do however lead to distinct spatio-temporal dynamics, selecting a speed of propagation of the instability and a state created in its wake. Curiously, in all ``generic'' examples that the authors are aware of, there is in fact a dichotomy of a ``positive'' front and a ``negative'' front, even in systems where there is no natural definition of positive and negative. As a consequence, only \emph{two} different states can be created in the wake of an invasion process. We shall explain in this introduction why it is reasonable to conjecture that such a dichotomy holds quite generally for pulled fronts, that is, for fronts whose dynamics are determined by the linear behavior in the leading edge, which in turn comes with just two possible directions in a leading eigenspace. It is this determinacy of front dynamics by the leading edge that in fact leads to the terminology of a pulled front, as opposed to pushed fronts whose speed is determined by the state in the wake.

Our main results, Propositions~\ref{p:Nfronts} and \ref{p:3fronts}, refute this conjecture of a more general dichotomy between positive and negative fronts determined by the leading edge linear dynamics by exhibiting an open (hence robust) class of systems that possess fronts whose speed is linearly determined, but which may leave { an arbitrary number of} different states in their wake depending on initial conditions. The example therefore demonstrates that dominant leading edge behavior, while responsible for the speed, does not select the state in the wake: pulled fronts are not just pulled.

\paragraph{The spreading dichotomy in a scalar example.} To be more precise, consider the simplest examples which arises in scalar equations of the type
\begin{equation}\label{e:nag}
u_t=u_{xx}+f(u),\qquad x\in\R,\ u\in\R,
\end{equation}
with a trivial unstable equilibrium at the origin, $f(0)=0$, and $f'(0)>0$, for instance  $f(u)=u(1-u)(1+u)$. Compactly supported, positive initial conditions $u_0(x)$ lead to solutions that converge to $u=1$ locally uniformly, or, more precisely,
\[
\lim_{t\to\infty} u(t,x-ct)=
\left\{\begin{array}{ll}
1,& |c|<c_\mathrm{lin}t,\\
0,& |c|>c_\mathrm{lin}t,
\end{array}\right.
\qquad c_\mathrm{lin}=2\sqrt{f'(0)}=2,
\]
while negative, compactly supported  initial conditions lead to convergence to $-1$ when $|c|<c_\mathrm{lin}$. More generally compactly supported initial conditions converge to either $+1$ or $-1$ along rays $x=ct$ for $0<|c|<c_\mathrm{lin}$. Regions where $u=+1$ and where $u=-1$ in the wake of the front are separated by kinks and anti-kinks $u(x)=\pm\tanh(x/\sqrt{2})$ \cite{PP22}.

\begin{figure}
\includegraphics[height=.8in]{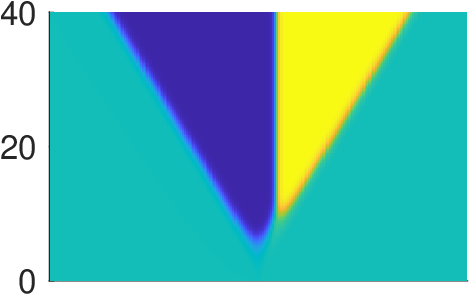}\hfill
\includegraphics[height=.8in]{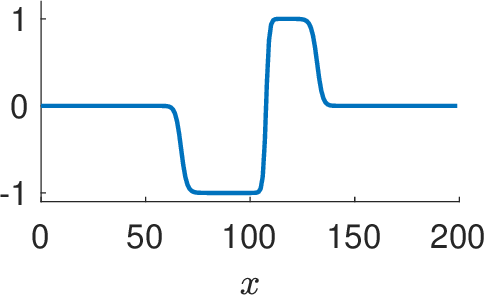}\hfill
\includegraphics[height=.8in]{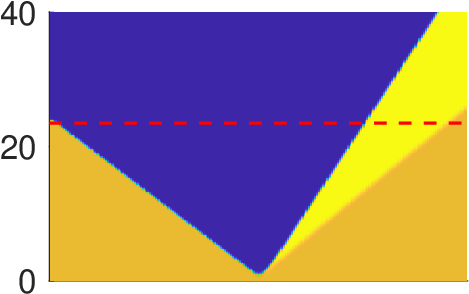}\hfill
\includegraphics[height=.8in]{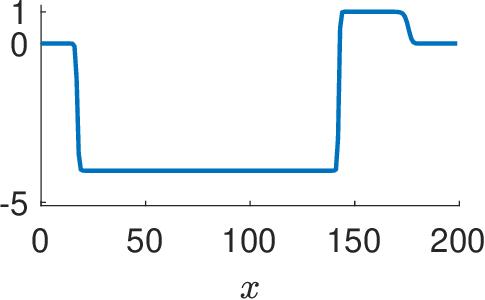}
\caption{Positive and negative fronts in the Nagumo equation \eqref{e:nag} originating from a sign-changing initial condition at the center of the domain:  balanced nonlinearity $a=1$ (left two panels) and imbalanced nonlinearity $a=.2$ (right two panels). Space-time plots and snapshot at $t=20$ in both case. Note the slightly different speeds of propagation in the imbalanced case due to a pushed front propagating to the left in the right two panels: the left front hits the boundary before the right front (illustrated by the horizontal red dashed line). Invasion fronts leave behind a stationary (left) or traveling (right) kink.}\label{f:ppselect}
\end{figure}

For weakly  asymmetric cubic nonlinearities, $f(u)=u(1-u)(u-a)$, with $-1<a<-1/2$, an equivalent result holds, that is, $u$ converges to $1$ or to $a$ along rays with speeds $|c|<2\sqrt{-a}$, with the exception of kinks now propagating with nonzero speed. Strong asymmetry $-1/2<a<0$ leads to different speeds of propagation: negative initial conditions spread with speed $|c|=2\sqrt{-a}$, positive initial conditions spread faster with speed $|c|=(1-2a)/\sqrt{2}$; see Figure~\ref{f:ppselect}.

A somewhat simpler description of the dynamics focuses on step-like initial data, $u=u_0$ for $x<0$, $u=0$ for $x>0$. For $u_0$ positive, solutions converge to a front leaving $u=1$ in its wake; for $u_0$ negative, the front leaves $u=-a$ in its wake.

\begin{figure}
\centering
\includegraphics[height=1.in]{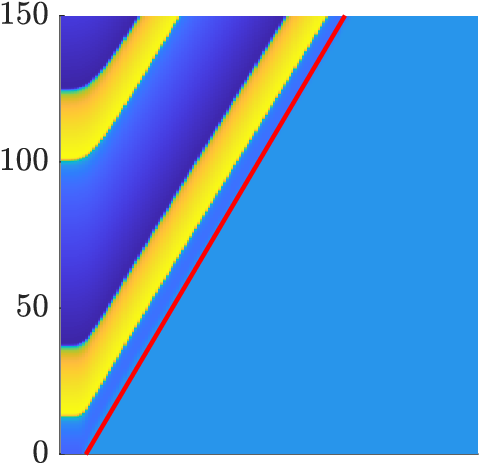}
\hspace*{0.07in}\includegraphics[height=0.90in]{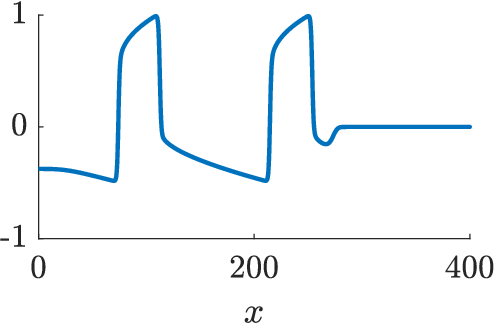}
\hspace*{.7in}
\includegraphics[height=1.in]{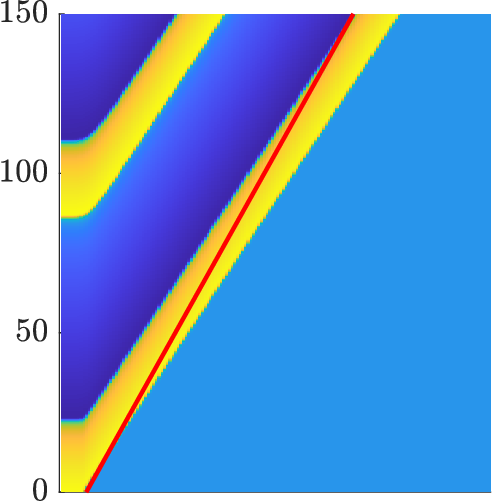}
\hspace*{0.07in}\includegraphics[height=0.90in]{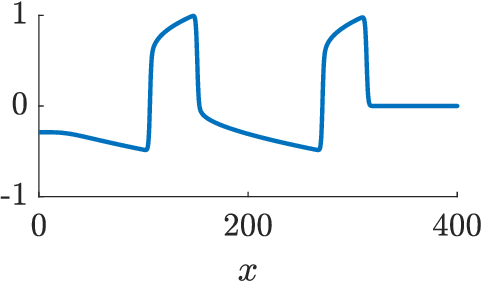}
\caption{Space-time plots  of the $u$-component of pattern-forming fronts in FHN \eqref{e:fhn}; parameters are $a=-0.2$, $b=\gamma=0$ showing a pulled front with negative leading edge (left) and a pushed front with positive leading edge (right).
Red lines indicates theoretical linear spreading speeds. Note the faster, pushed, speed in the right panel.
Also shown are profiles of solution $u(t,x)$ at times $t=150$). 
%
}
\label{f:fhn}
\end{figure}

\paragraph{Dichotomies in systems.} Somewhat curiously, this dichotomy between two possible modes of invasion is observed in many \emph{systems} of reaction-diffusion equations, and one may, in fact,  even expect such a dichotomy to hold more universally. We are thinking here for instance of systems     
\begin{equation}\label{e:rd}
U_t=DU_{xx}+F(U),\qquad x\in\R,\ U\in\R^N,\ D\in\R^{N\times N}>0
\end{equation}
again with trivial unstable equilibrium at the origin $F(0)=0$, $\mathrm{spec}\,F'(0)\cap\{z\in\C|\Re z>0\}\neq 0$. One often finds spreading of step-like initial data at a distinct speed, in particular when the system obeys a comparison principle and initial data is positive (or negative); see for instance \cite{llw} for cooperative systems with a comparison principle and \cite{cartersch,avery2023stabilitycoherentpatternformation,PushedFHN} for a dichotomy in the FitzHugh-Nagumo equation,
\begin{equation}\label{e:fhn}
\begin{aligned}
u_t&=u_{xx}+u(1-u)(u-a) -v,\nonumber\\
v_t&=\eps(u-\gamma v + b),\nonumber
\end{aligned}
\end{equation}
which does not possess a comparison principle; see Figure~\ref{f:fhn}.

%

\paragraph{The linear origin of the dichotomy.} In order to determine a speed and a state $U_-$  selected in the wake, one starts with the dynamics in the leading edge of the invasion process, which are governed by the linearized equation
\begin{equation}\label{e:rdlin}
U_t=DU_{xx}+F'(0)U.
\end{equation}
Spreading speeds are determined by marginal pointwise stability, that is, the spreading speed is the largest speed in which one observes pointwise instability; see \cite{HolzerScheelPointwiseGrowth} for background. Pointwise growth and decay in a comoving frame are (generically) determined by the presence of pinched double roots $(\lambda_\mathrm{dr},\nu_\mathrm{dr})$ on the imaginary axis, that is solutions to
\begin{equation}\label{e:dr}
d_c(\lambda,\nu)=0,\quad \partial_\nu d_c(\lambda,\nu)=0,
\end{equation}
where
\[
d(\lambda,\nu)=\mathrm{det}\,(D\nu^2+F'(0)-\lambda),\qquad d_c(\lambda,\nu)=d(\lambda-c\nu,\nu),
\]
together with the pinching condition which encodes that the two roots $\nu_\pm(\lambda)$ which collide at $(\lambda_\mathrm{dr}, \nu_\mathrm{dr})$ originate from opposite directions as $\mathrm{Re} \, \lambda \to \infty$. More precisely, the pinching condition requires that two continuous curves of roots $\nu=\nu_\pm(\lambda(\tau))$, $\tau\in[0,\infty)$ to $d_c(\lambda,\nu))=0$, $\nu_\pm(\lambda)$ on a path $\lambda(0)=\lambda_\mathrm{dr}$,  and $\Re\lambda(\tau)$ increasing, $\Re\lambda(\tau)\to\infty$ for $\tau\to\infty$, with
$\nu_\pm(\lambda(0))=\nu_\mathrm{dr}$, satisfy
\[
\Re\nu_\pm(\lambda(\tau))\to\pm\infty, \qquad \text{for }\tau\to\infty.
\]
Generically $(\partial_\lambda d_c)\cdot (\partial_{\nu}{\nu}d_c)\neq 0$ and there is a unique eigenvector $U_\mathrm{dr}$ and generalized direction $U^1_\mathrm{dr}$ associated with the pointwise growth, that is, for $\lambda=\lambda_\mathrm{dr}, \nu=\nu_\mathrm{dr}$,
\[
\begin{aligned}
\lambda U_\mathrm{dr} &=D\nu^2 U_\mathrm{dr}+c\nu{U}_\mathrm{dr}+F'(0)U_\mathrm{dr},\\
\lambda U_\mathrm{dr}^1 &=D\nu^2 U_\mathrm{dr}^1+c\nu{U}_\mathrm{dr}^1+F'(0)U_\mathrm{dr}^1+(2\nu+c)U_\mathrm{dr}.
\end{aligned}
\]
We focus on the situation where $\lambda_\mathrm{dr}=0$ at $c=c_\mathrm{lin}$, and $\Re\nu_\mathrm{dr}<0$, so that the behavior in the leading edge of the instability is non-oscillatory. The vector $U_\mathrm{dr}$ then gives the direction in $\R^N$ of the instability, complementary directions decay in the leading edge. This suggests that the linearized instability points in precisely two directions, $+U_\mathrm{dr}$ and $-U_\mathrm{dr}$.

Invasion processes that are determined by these linear dynamics, since for instance nonlinear effects do not enhance the speed of propagation, are commonly referred to as pulled and the associated fronts, which propagate with the associated linear speed of propagation and are marginally stable precisely due to the marginal stability in the leading edge that we just analyzed are referred to as pulled fronts. The idea behind this terminology is that the dynamics in the leading edge ``pull'' the instability into the unstable state.
This is the linear origin of the conjecture of a dichotomy of fronts, determined by a sign in the leading edge corresponding to ``positive'' perturbations $+U_\mathrm{dr}$ and  ``negative'' perturbations $-U_\mathrm{dr}$, respectively.

\paragraph{A trivial (degenerate) counterexample.} There are examples where such a simple dichotomy can obviously not be true, since the leading edge does not possess a unique eigenvector $U_\mathrm{dr}$, that is, the double root is degenerate as a solution to \eqref{e:dr}. A simple such example is the complex Ginzburg-Landau equation,
\begin{equation}\label{e:cgl}
A_t=A_{xx}+A-A|A|^2,
\end{equation}
which we write here in real variables as
\begin{equation}\label{e:rcgl}
\begin{aligned}
u_t&=u_{xx}+u(1-u^2-v^2),\\
v_t&=v_{xx}+v(1-u^2-v^2).
\end{aligned}
\end{equation}
The linearization decouples into two identical copies of the exponentially growing diffusion equation $u_t=u_{xx}+u$ with spreading speed $2$. In fact propagation at any complex angle is possible and equally favored, that is, one can observe fronts $A(t,x)=\rme^{\rmi\varphi} u(x-2t)$ for any $\varphi\in[0,2\pi)$, $u(\xi)\in\R$, $u(\xi)\to -1$ for $\xi\to -\infty$ and $u(\xi)\to 0$ for $\xi\to\infty$ \cite{EBERT200413,as5,EW1994}.

The Ginzburg-Landau equation arises as a normal form or modulation equation for pattern-forming systems near a Turing instability \cite{crosshohenberg,mielke}, the simplest example being the Swift-Hohenberg equation, where the instability is in fact oscillatory, leading to the rhythmic nucleation of patterns in the leading edge. The modulation is however averaged out, leaving a double zero eigenvalue at the origin. From our perspective here, this situation is degenerate and we wonder if perturbations that break the degeneracy would result in the selection of precisely two selected fronts, a positive and negative front associated with the eigenvector to the pinched double root system.

\paragraph{Counterexample: more than two fronts.} To rule out examples such as the complex Ginzburg-Landau equation \eqref{e:cgl} which are in some sense special, or degenerate, we insist that counter examples should be robust. Our main result states that { given a natural number $N$,} there is in fact an open class of systems that exhibit {$N$ selected fronts}. To fix the setting, consider therefore
\begin{equation}
\begin{aligned}
u_t&=d_u u_{xx} + f_u(u,v),\\
v_t&=d_v v_{xx} + f_v(u,v),
\end{aligned}
\end{equation}
where we think of $d_{u/v}\in\R^2$ and $f_{u/v}\in C^2_\mathrm{loc}(\R^2,\R)$ as parameters. To more easily state our result, we fix the origin as an equilibrium and define $\mathcal{F}=\{f_{u}(0,0)=f_{v}(0,0)=0\}\subset C^2_\mathrm{loc}(\R^2,\R^2)\}$.

\begin{proposition}[Multiple fronts]\label{p:3fronts}
For any $N\in\N$, there exists an open region $\Omega\subset \R^2\times \mathcal{F}$, so that  for all $(d_u,d_v,f_u,f_v) \in \Omega$, \eqref{e:cglf} possesses $N$ distinct selected, pulled fronts $(u,v)_{\mathrm{f},j}(x-c_\mathrm{lin}t)$, { $j=1,2, ..., N$}; see Definition~\ref{d:sel}, below for a precise characterization of selected, pulled fronts.
Moreover, $\lim_{\xi\to-\infty} (u,v)_{\mathrm{f},j}(\xi)=(u,v)_{*,j}$ are all different, that is, the fronts in fact select different states in their wake.
\end{proposition}
The proposition refers to the selected pulled front, a notion that we define next. We define what we mean by selection, first, and then list specific criteria that guarantee selection. Those criteria are commonly associated with the notion of a pulled front. Pushed fronts are characterized by a different set of conditions, encoding marginal stability due to a localized mode at the front interface, which also guarantees selection. 
%
%
%
%

\begin{definition}[Selected fronts]\label{d:sel}
We say a front $U_*$ with speed $c_*$ is selected if there exists a normed space $X$ and for any $\eps>0$ an open class $\mathcal{U}_\eps\subset X$ of initial conditions, containing functions with support in $\R_-$, so that solutions $U$ stay close to the family of translates, $\|U(t,\cdot { + c_* t})-U_*(\cdot-h(t))\|_{L^\infty}<\eps$, for $t$ sufficiently large, with a phase shift $h(t)=\rmo(t)$ as $t\to\infty$.
\end{definition}
The topology $X$ typically enforces exponential and algebraic decay near $x=+\infty$. The definition was introduced in \cite{as1} to identify speeds and front profiles that are relevant for spreading from compactly supported data. The shift $h(t)$ is necessary due to a slow drift along front profiles, where typically $h(t)\sim -\log t$; see \cite{Bramson2,Lau,EBERT200413,Comparison1}.

Selection in the sense of Definition~\ref{d:sel} was established in \cite{as1} for scalar parabolic equations and in \cite{avery2} for parabolic systems relying on robust properties of front profiles that we shall make explicit, next.

\begin{definition}[Marginal stability and pulled fronts]\label{d:stab}
We say a front $U_*(x-c_\mathrm{lin}t)$, $U_*(\xi)\to 0$ for $\xi\to+\infty$, $U_*(\xi)\to U_-$ for $\xi\to-\infty$, with speed $c_\mathrm{lin}$ satisfies the pulled front selection criterion if the linearization $\mathcal{L}$ is marginally stable in an exponentially weighted norm, with marginally stable spectrum only due to a pinched double root at $\lambda=0$ in the leading edge, that is,
\begin{itemize}
\item $U_-$ is linearly exponentially attracting;
\item $(0,-\eta)$ with $\eta>0$ is a pinched double root for the linearization at $U=0$;
\item $\mathrm{spec}_\eta\,(\mathcal{L})\cap \{\Re\lambda\geq0\}=\{0\}$, where the norm is $\|U\|_\eta=\|(1+\rme^{\eta\cdot })U(\cdot)\|_{L^\infty}$;
\item the kernel of $\mathcal{L}$ in the norm does not contain functions with $\|\cdot\|_\eta<\infty$, but $U_*(\xi)\rme^{\eta\xi}/\xi\to a_+\neq 0$.
\end{itemize}
\end{definition}

A similar definition would characterize pushed fronts as marginally stable due to only a simple eigenvalue, rather than essential spectrum, at the origin. The results in \cite{as1,avery2} then give the following.
\begin{theorem}[Marginal stability $\Longrightarrow$ selection]\label{t:ms}
Fronts satisfying the selection criteria of Definition~\ref{d:stab} are selected in the sense of Definition~\ref{d:sel}.
\end{theorem}

As a corollary, one also finds robustness \cite{as1,avery2}.

\begin{theorem}[Robustness of marginal stability]\label{t:rob}
Fronts satisfying the selection criteria of Definition~\ref{d:stab} occur for open sets of reaction-diffusion equations. 
\end{theorem}

As a consequence, in order to prove Proposition \ref{p:3fronts}, it is sufficient to exhibit an example that is marginally stable in the sense of Definition \ref{d:stab}. We give two such examples. In the first example, which may appear somewhat artificial, we verify all assumptions analytically. In the second example, we show how the first example underlies a more general mechanism for multiple front selection in a more realistic setting.

\paragraph{Example 1: Multistability in scalar equations and multiple fronts.}

Our first example is a simple scalar equation
\begin{equation}
u_t=u_{xx}+f(u),
\end{equation}
with a carefully constructed $f$.

\begin{proposition}[$N$ fronts]\label{p:Nfronts}
For each $N\in\N$ There exists an open class of nonlinearities $f\in \{C^2_\mathrm{loc},f(0)=0\}$ such that there exist $N$ selected invasion fronts leaving behind states $u_j$, $j=1,\ldots,N$, which are all distinct $u_0=0<u_1<u_2<\ldots<u_N$.
\end{proposition}
The states $u_j$ are all stable, $f'(u_j)<0$, while $f'(0)>0$ leading to a linear spreading speed $c_\mathrm{lin}=2\sqrt{f'(0)}$. The nonlinearity is constructed as a perturbation of a nonlinearity which allows for monotone fronts all traveling with speed $c_\mathrm{lin}$, of the form
\[
\lim_{x\to +\infty} u_{\ell,*}(x)=u_{\ell-1}, \lim_{x\to -\infty} u_{\ell,*}(x)=u_\ell,\qquad \text{ for } \ell=1,\ldots, N.
\]
The perturbation gives fronts connecting $u_\ell$ to $0$ by concatenating fronts $u_{\ell,*}$,\ldots,$u_{1,*}$. 

Situations of this type were more generally analyzed for instance in \cite{gilettimatano} and referred to as propagating terraces. In out particular interest, we require all fronts, that is, all ``steps'' in the terrace to propagate at the same speed.

\paragraph{Example 2: Skew-coupled staged invasion.}
We next turn to an example where the plateaus in between fronts are unstable states. These turn out to be more common in applications, but occur in systems of equations, only. We consider the skew-coupled reaction-diffusion system,
\begin{equation}\label{e:sk}
\begin{aligned}
    u_t&=u_{xx}+u-u^3,\\
    v_t&=v_{xx}+(2u^2-1+\mu)(v-v^3).
\end{aligned}
\end{equation}
Setting $v=0$, we find a simple invasion problem with unstable state $u=0$. There is a selected front $u_*(\xi)\to 1$ for $\xi\to -\infty$,  $u_*(\xi)\to 0$ for $\xi\to \infty$, where $\xi=x-2t$.  There also is, by symmetry, a second selected front $-u_*(\xi)$. Restricted to the first equation, only, both these fronts are marginally stable in the sense of Definition \eqref{d:stab}. The fronts are of course also solutions to the full system, since $v=0$ is simply invariant under the time dynamics. Linearizing at the front in the full equation however yields an instability, so that as solutions to the full system, the fronts do not satisfy the conditions in Definition \ref{d:stab}. To see this, note that the linearization at the front in a comoving frame (here, denoted by $x$ again) is block-triangular, and simply inspect the linearization at the front in the $v$-equation,
\[
v_t=v_{xx}+2v_x+(2u_*^2-1+\mu)v=:\mathcal{L}_v v
\]
Since $u_*^2\to 1$ for $x\to -\infty$, we have that the essential spectrum of $\mathcal{L}_v$ contains  the spectrum of $\mathcal{L}_-=\partial_{xx}+2\partial_x+1+\mu$, with maximal real part $1+\mu$. In an optimal exponential weight $\|v\|_{L^2_\eta}:=\|w\|_{L^2}$, with $w(x)=v(x)\rme^{\eta x}$, $\eta=1$ for maximal stability,  the spectrum of $\mathcal{L}_-$ is $(-\infty, \mu]$, that is, increasing $\mu$ past zero will trigger an instability in any exponentially weighted space.

Somewhat equivalently, we can fix $u=1$ (or $u=-1$), the state left in the wake of the $u$-invasion front and consider the $v$-equation for this fixed $u$-value, 
\[
v_t=v_{xx}+(1+\mu)(v-v^3).
\]
Here, $v=0$ is unstable and the instability spreads with speed $c=2\sqrt{1+\mu}$, mediated by fronts $\pm v_*(x-ct;c)$,  $v_*(\xi;c)\to 1$ for $\xi\to-\infty$, $v_*(\xi)\to 0$ for $\xi\to +\infty$, $\xi=x-ct$. In particular, the instability spreads faster than the $u$-instability for $\mu>0$, and slower for $\mu<0$. It is then reasonable to expect that for $\mu\gtrsim 0$, there is a front consisting of a concatenation of the $u_*$, ramping up $u$ from $0$ to $\pm 1$ followed by a $v$-front, ramping up $v$ from $0$ to $\pm 1$. These fronts are in fact easily observable numerically. We will show below that these four fronts, leaving $(u,v)=(\pm 1,\pm 1)$ with any combination of signs possible in their wake, in fact exist and satisfy Definition \ref{d:stab}.

\begin{proposition}\label{p:4fronts}
For $\mu\gtrsim 0$, there exist four invasion fronts $(\pm u_*(\xi),\pm v_*(\xi))$ leaving the states $(u,v)=(\pm 1,\pm 1)$ in their wake. All four fronts are marginally stable in the sense of Definition \ref{d:stab}.
\end{proposition}
We emphasize that the robustness result, Theorem \ref{t:rob}, ensures that there are also 4 fronts for perturbations of the specific system \eqref{e:sk} that include dependence on $v$ in the $u$-equation or additive $u$-dependent coupling terms in the $v$-equation, thus destroying the skew-coupling structure and the invariance of subspaces $\{u=0\}$ or $\{v=0\}$.

This example can easily be combined with the scalar example, replacing the kinetics in the $u$-equation, to obtain any number of invasion fronts, {thus establishing our main result, Proposition \ref{p:3fronts}. We nonetheless continue exploring different aspects of non-uniqueness of selected fronts with the following examples, which may be more prevalent in application.}

\paragraph{Example 3: Forced complex Ginzburg-Landau.}
Our next example builds on the idea of secondary invasion fronts creating alternative states. The system we consider is a perturbation of \eqref{e:cgl} that arises when studying pattern formation in spatially prepatterned systems, with prepatterning at one half and one third of the wavelength selected by the instability; see \cite{PhysRevLett.51.786,PhysRevLett.87.238301,PhysRevLett.56.724}.
Specifically, we consider
\begin{equation}\label{e:cglf}
    A_t=A_{xx}+A-A|A|^2 +\alpha \bar{A}+ \beta \bar{A}^2.
\end{equation}
Rotating $A\mapsto \rme^{\rmi\varphi}A$, we may assume $\alpha>0$. While there are many interesting questions associated to the relative complex phase of $\alpha$ and $\beta$, we focus here on $\beta>0$. This leads to the real variable equation
\begin{equation}\label{e:rcglf}
\begin{aligned}
u_t&=u_{xx}+(1+\alpha) u +\beta (u^2-v^2)-u(u^2+v^2),\\
v_t&=v_{xx}+(1-\alpha) v -2\beta uv-v(u^2+v^2).
\end{aligned}
\end{equation}
For $\alpha\neq 0$, the linear spreading speed is $2\sqrt{1+\alpha}$. In the equation linearized at the origin,  with $\alpha>0$, the $u$-component spreads at this speed while the $v$-component decays exponentially pointwise in the frame moving with this speed.

\begin{figure}
\centering\includegraphics[width=0.7\textwidth]{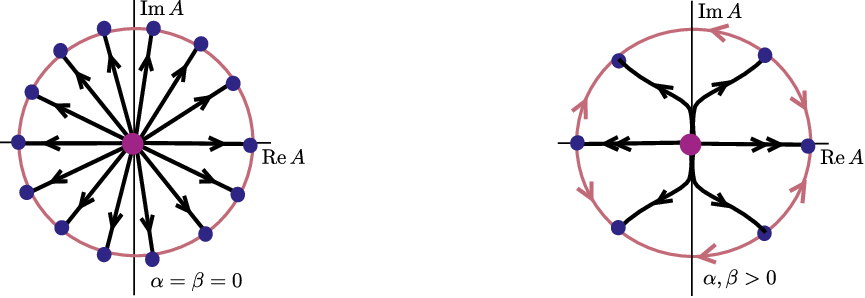}
\caption{Schematic phase portrait of \eqref{e:cglkinf} for
$\alpha=\beta=0$ (left) and $\alpha,\beta \gtrsim 0$ (right). Note that equilibria on the invariant (almost) circle are alternatingly stable and unstable.} \label{f:cglkin}
\end{figure}

The $x$-independent dynamics,
\begin{equation}\label{e:cglkinf}
    A_t=A-A|A|^2 +\alpha \bar{A}+ \beta \bar{A}^2\in\C,
\end{equation}
possess $6$ equilibria $A_{j}\sim \rme^{2\pi\rmi j/6}$, $j=0,\ldots,5$  near the circle $|A|=1$. The equilibria $A_{0,2,4}$ are stable, the other equilibria are saddles; see Figure~\ref{f:cglkin}.
We claim that all three stable nontrivial states can be selected in the wake of an invasion front.
We only provide numerical evidence and heuristics based on the conceptual gluing procedure in Example 1 and therefore state this  fact as a result rather than a proposition.

\begin{result}\label{res:3frontscgl}
For $\beta\in  (\beta_\mathrm{c},\beta_\mathrm{p})$, with $\beta_\mathrm{c}\leq \frac{1+3\alpha}{\sqrt{6}}$ and $\beta_\mathrm{p}=\sqrt{\frac{1+a}{2}}$, $0<\alpha<1/3$, there exist three pulled fronts satisfying the assumptions of Definition \ref{d:stab}, leaving the stable states $A_j$, $j=0,2,4$ in their wake.
\end{result}
The front leaving $A_0$ in its wake is simply the positive  real front. The other two fronts arise from a secondary invasion: a primary negative real front leaves $A_3<0$ in its wake. This state is unstable against perturbations in the imaginary direction. Postivise and negative imaginary perturbations spread at a speed faster than the primary speed for $\beta>\frac{1+3\alpha}{\sqrt{6}}$. The upper boundary $\beta=\sqrt{\frac{1+a}{2}}$ marks a transition of the primary front to a pushed front.

\begin{figure}
    \centering\includegraphics[height=1.in]{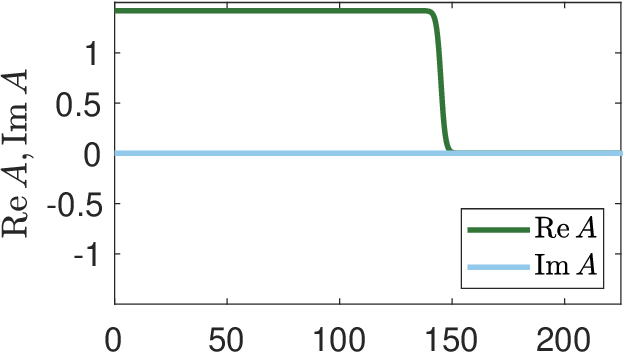}\qquad
    \includegraphics[height=1.in]{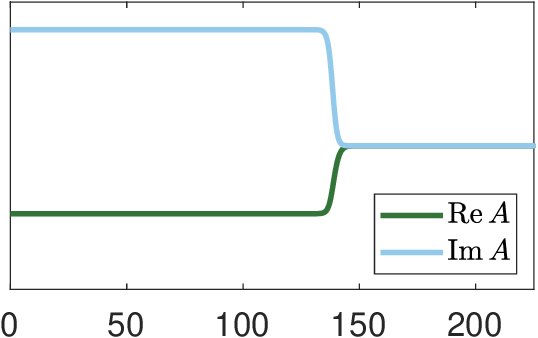}\qquad
    \includegraphics[height=1.in]{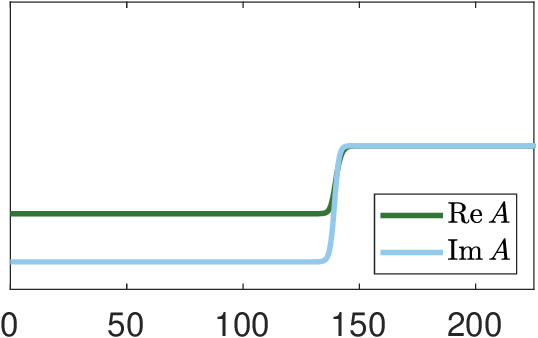}\\[0.2in]
   \includegraphics[height=1.in]{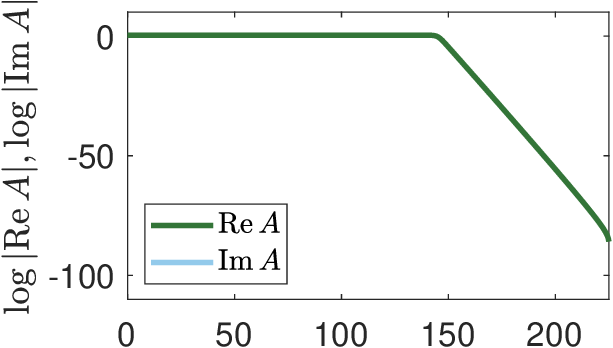}\qquad
    \includegraphics[height=1.in]{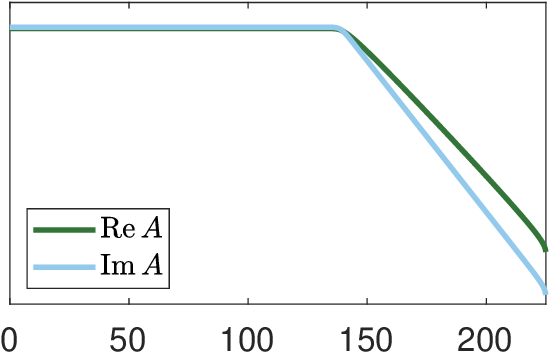}\qquad
    \includegraphics[height=1.in]{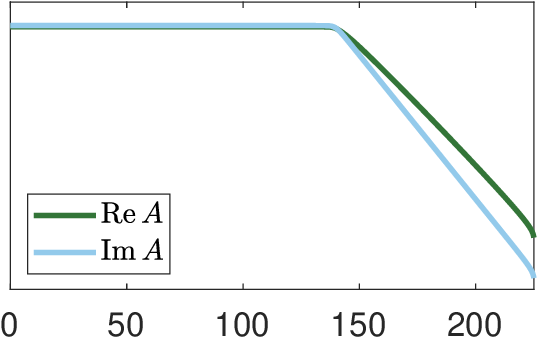}
    \caption{Solution to \eqref{e:cglf} with $\beta=0.5$, $\alpha=0.02$, at time $T=1500$, in a frame moving with the linear spreading speed $c_*=2\sqrt{1+\alpha}$, and initial condition $A=\rme^{\rmi \ell\pi/3}\chi_{<0}$, $\ell=0,1.2,4$ from left to right. All three fronts are selected and pulled.}\label{f:threefronts}
\end{figure}

All three fronts can be observed starting from step-function initial data with the state $A=\rme^{2\pi \rmi j/6}$, $j=0,2,4$ in $x<0$. Figure~\ref{f:threefronts} shows the observations in a comoving frame with speed $c=2\sqrt{1+\alpha}$ confirming our predictions.

\begin{figure}
        \centering\includegraphics[width=.26\textwidth]{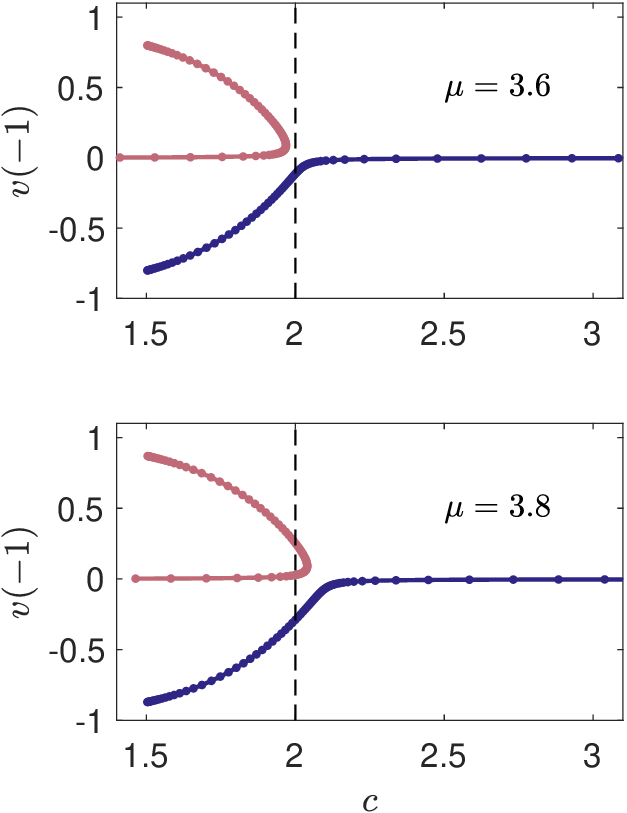}\includegraphics[width=.26\textwidth]{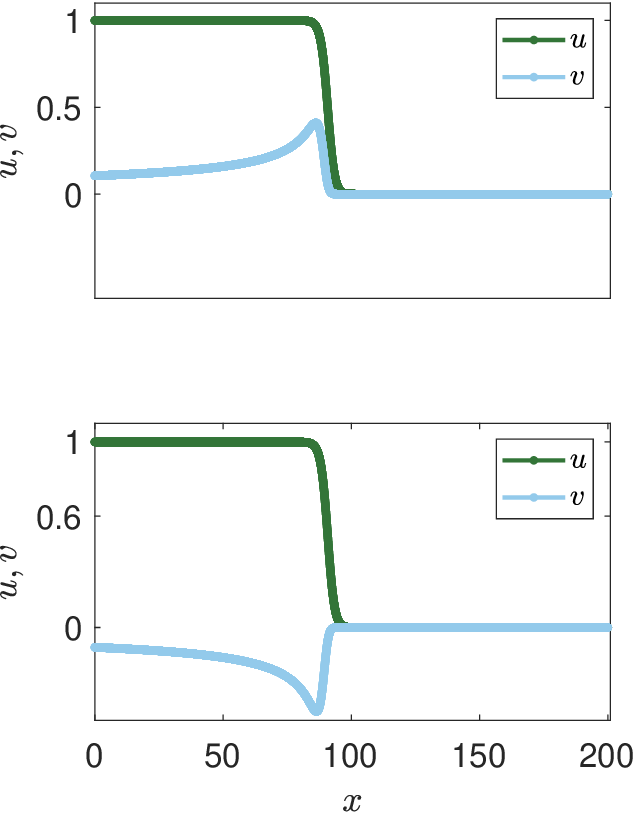}\includegraphics[width=.2\textwidth]{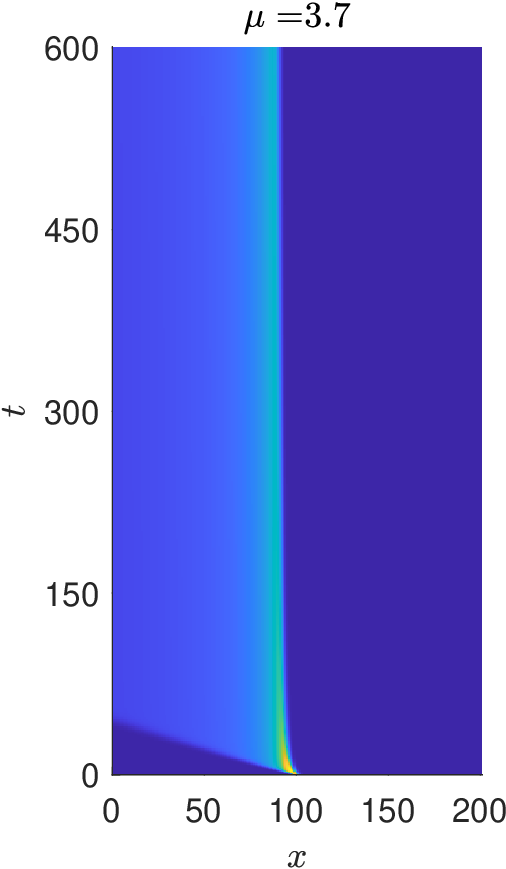}\includegraphics[width=.2\textwidth]{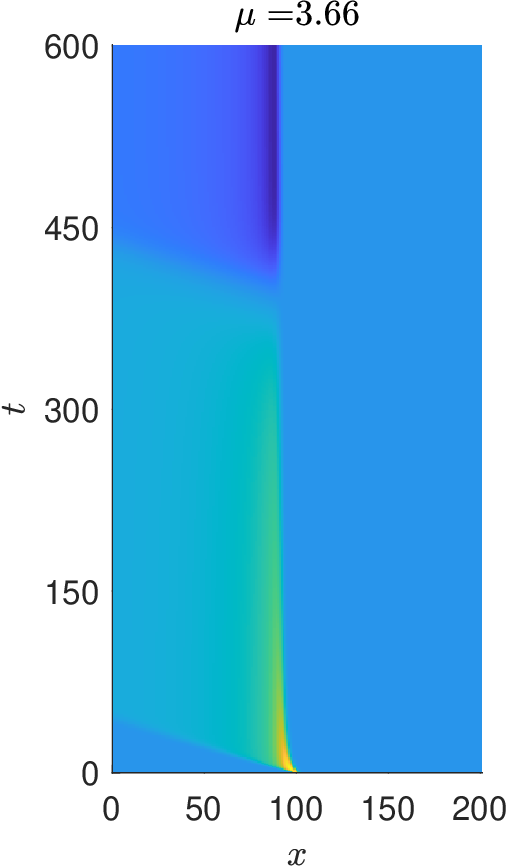}
        \caption{A saddle-node of invasion fronts occurs in \eqref{e:losstrans} with $\delta=10^{-3}$. The saddle occurs when the speed is varied (left) for a value $c_\mathrm{sn}(\mu)<2$ for $\mu=3.6$ and $c_\mathrm{sn}(\mu)>2$ for $\mu=3.8$ (left); front profiles for $u$ and $v$ (center, $\mu=3.7$) and space-time plots of $v$ (right) from direct simulations of the invasion process in a frame traveling with speed 2 illustrate that, fixing $c=2$, we see a saddle-node of selected fronts: the positive fronts persists for a long time as the remainder of the equilibrium that has in fact disappeared in a saddle-node ($\mu=3.66$). } \label{f:frontsn}
\end{figure}

\paragraph{Example 4: Interface bifurcations in skew-coupled systems.}
Our last example shows non-uniqueness of fronts leaving behind the same state, as opposed to all previous examples where different fronts left different states in their wake. The multiplicity arises through eigenmodes that cause bifurcations localized in the front interface, without leading to pushed fronts. We consider
\begin{equation}\label{e:losstrans}
    u_t=u_{xx}+u-u^3,\qquad v_t=v_{xx}+\mu (u-u^3)v-v^3 + \delta (u-u^3).
\end{equation}
Here, the $u$-equation clearly possesses front solutions connecting $u=1$ to $u=0$ for all speeds. The $v$-equation, at $\delta=0$, possesses the trivial solution $v\equiv 0$ with a linearization of the form $\partial_{xx}+c\partial_x +V(x)$ with localized, positive potential $V(x)=\mu(u-u^3)$, where $u$ is the front with speed $c$. It turns out that for large speeds $c$, this localized instability is destroyed by the advection term, while for small speeds $c$ it is present and leads to a pitchfork bifurcation and the emergence of a localized pattern in the $v$-equation, emerging near the interface of the $u$-front. This localized pattern can be either positive or negative due to the reflection symmetry in the $v$-equation. The critical speed $c_\mathrm{bif}$ depends on $\mu$, and coincides with the spreading speed for roughly $\mu=3.6$. Adding a small perturbation $\delta\neq 0$ breaks the pitchfork into two branches: for negative $\delta$, a localized negative bump in the $v$-equation exists for all values of $c$. A solution with a positive bump in $c$ exists for $c<c_\mathrm{bif}(\mu)$, where it undergoes a saddle-node bifurcation. For $\delta$ small, there is a $\mu_\mathrm{bif}\sim 3.6$ so that $c_\mathrm{bif}(\mu_\mathrm{bif})=c_\mathrm{lin}=2$; see Fig. \ref{f:frontsn}.

In summary, the example above allows for the emergence of a nonlinear front with speed $c_\mathrm{lin}$ through a saddle-node bifurcation, which cannot be found through a homotopy in $c$ starting with $c\gg 1$. We will return to this example when discussing marginal stability, below.

\section{Proof of Proposition \ref{p:Nfronts} --- $N$ fronts in scalar PDE invasion}\label{s:ex0}
\begin{figure}
\centering
\includegraphics[width=0.9\textwidth]{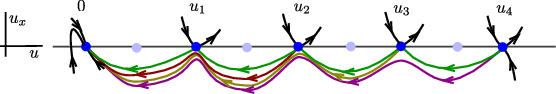}
\caption{Schematic picture of the equilibria $u_j$ (blue), heteroclinics in the unperturbed equation (green) and perturbed chains of heteroclinics. The perturbed heteroclinics do not lie in the strong stable manifold of the origin and are monotone, $u_x<0$. }\label{f:scal_ex}
\end{figure}

We construct $f$ with the desired properties. We start with a function $f$ such that $f(\ell)=0$, $\ell=0,1,2,3,\ldots N$, and $f'(0)>0$, $f'(\ell)<0$, $\ell=1,\ldots,N$.
Let $c=c_\mathrm{lin}=2\sqrt{f'(0)}$ be the linear spreading speed. The goal is to establish a phase portrait for   the ODE
\[
u_x=v,\qquad v_x=-c_\mathrm{lin}v-f(u).
\]
as depicted in Figure \ref{f:scal_ex}. We start with the unperturbed phase portrait, which shall  exhibit heteroclinic connections $u_{\ell+1}\to u_\ell$ for all $\ell=0,\ldots, N-1$, shown in green in Fig. \ref{f:scal_ex}. We therefore modify $f$ in closed intervals $J_\ell\subset (\ell-1,\ell)$, $\ell=1,\ldots,N$ so that there indeed  exist monotone heteroclinic connections $u_{\ell,*}$ connecting $u_\ell$ to $u_{\ell-1}$, at speed $c_\mathrm{lin}$. Moreover, we may assume that the heteroclinic $u_{1,*}$ is a non-degenerate pulled front, that is, the heteroclinic connection is not contained in the (codimension-one) strong stable manifold, consisting of trajectories with pure exponential asyptotics $\rme^{-\sqrt{f'(0)}x}$. We now successively modify the nonlinearity in the intervals $(\ell,\ell+1)$, $\ell=1,\ldots,N-1$, starting with $\ell=1$. We change the nonlinearity in such a way that the heteroclinic connection breaks and the unstable manifold of $u_2$ intersects the line $(u_1,0)$ in $u-u_x$-phase space at some small negative $(-\eps,0)$. Following this manifold further, it stays close to the heteroclinic connection from $u_1$ to $0$ and therefore is also contained in the stable manifold of the origin, excluding the strong stable manifold. In other words, upon perturbation, the unstable manifold of $u_2$ just misses the stable manifold of $u_1$. Continuity of the transition near the equilibrium $u_1$, which can for instance be seen in coordinates that linearize the dynamics near $u_1$ using the Grobman-Hartman theorem, then gives that the unstable manifold of $u_2$ is close to the unstable manifold of $u_1$ when it intersects a line $u=u_1-\eps$, just to the left of $u_1$. Closeness to the unstable manifold of $u_1$ when entering a neighborhood of the origin then also guarantees that the unstable manifold of $u_2$ also misses the strong stable manifold of the origin with pure exponential asymptotics. 
This construction therefore gives a pulled invasion front leaving behind $u_2$ in the wake, in addition to the front leaving behind $u_1$. We now ignore the front leaving behind $u_1$ and continue our construction concatenating the front connecting $u_3$ and $u_2$ with the newly found invasion front connecting $u_2$ and the origin, repeating exactly the above construction. Again, modifying the nonlinearity in $J_3$, we can find a front connecting $u_3$ and $0$ directly. Inductively, we can thus find $N$ pulled invasion fronts.

\section{Proof of Proposition \ref{p:4fronts} --- four fronts in staged invasion}\label{s:ex}

\paragraph{A short proof based on comparison principles.}

We consider the $v$-equation in a frame moving with speed $2$ with $\mu>0$ and $u=u_*$ fixed as the invasion front in the $u$-equation with, say, $u_*(0)=1/2$, 
\begin{equation}\label{e:vskew}
v_t=v_{xx}+2 v_x+ (2u_*^2-1+\mu)(v-v^3).
\end{equation}
We construct sub and super-solutions. Fix $\mu>0$ and choose any $0<\eps<\mu$.  Then there exists a $x_0\ll -1$ such that $2u_*(x_0)^2-1+\mu>1+\eps$.  Monotonicity of $u^*$ then implies that $2u_*(\xi)^2-1+\mu>1+\eps$  for all $x<x_0$.  Next, by essentially replacing the linear term in (\ref{e:vskew}) by $1$, let $\phi(\xi,c)$ denote the invasion front with speed $c$ which satisfies the equation
\[ 0=\phi''+c\phi'+\phi-\phi^3,\]
unique up to translation in $\xi$. When $c=2-\eps<2$ this front is oscillatory.  To define $\phi(\xi,c)$ uniquely we select the translate of the front that satisfies
\[ \phi(x_0,c)=0, \quad \phi(\xi,c)>0 \ \text{for all} \ \xi<x_0. \]
We then define the candidate sub-solution 
\[ \underline{v}_{\epsilon}(x,t)=\left\{ \begin{array}{cc} \phi(x+\eps t,2-\eps) & x+\eps t\leq x_0 \\ 0 & x+\eps t> x_0 \end{array}\right.,  \]
Define the operator
\[ N(v)=v_t-v_{xx}-2v_x-(2u^2_*-1+\mu)(v-v^3). \]
Then for $x+\eps t<x_0$
\begin{equation}
N\left(\underline{v}_\eps(x,t)\right) = (\phi-\phi^3)\left(1-(2u^2_*-1+\mu)\right) < -\eps (\phi-\phi^3) <0.
\end{equation} 
This establishes that $\underline{v}_\eps(x,t)$ is a sub-solution when $x+\eps t<x_0$ and, by extension, for any $(x,t)\in\mathbb{R}\times\mathbb{R}^+$.

For the super-solution we first consider the linear equation 
\[ w_t=w_{xx}+2w_x.  \]
For any $\nu<0$, this equation has solutions of the form $e^{\lambda t+\nu x}$ provided that 
\[ \lambda=\nu^2+2\nu. \]
These exponential solutions propagate with speed
\[ c_{\mathrm{env}}(\nu)=-\nu-2. \]
Take $\nu=-2-\eps$ so that $c_{\mathrm{env}}(-2-\eps)=\eps$ and the linear solution propagates to the right (in the co-moving frame adapted to the $u$ front).

Define
\[ \overline{v}_\eps(t,x)=\left\{ \begin{array}{cc} 1 & x\leq \eps t +\frac{1}{2+\eps} \mathrm{ln}(M) \\
M e^{\nu (x-c_{\mathrm{env}}t)} & x> \eps t +\frac{1}{2+\eps} \mathrm{ln}(M)\end{array}\right. ,  \]
for a constant $M$ to be defined below.  We then compute 
\[ N\left(M e^{\nu (x-c_{\mathrm{env}}t)}\right)= (1-\mu-2u^2_*)\left( Me^{\nu (x-c_{\mathrm{env}}t)}-\left(Me^{\nu (x-c_{\mathrm{env}}t)}\right)^3\right)    \]
Since  $u_*\to 0$ as $x\to\infty$ we can find an $x_*$ such that $1-\mu-2u^2_*>0$ for all $x>x_*$.    Finally, by taking $M$ sufficiently large so that $x_*<\frac{1}{2+\eps}\mathrm{ln}(M)$  we obtain that $\overline{v}_\eps$ is a super-solution for any $(x,t)\in\mathbb{R}\times\mathbb{R}^+$.

For any small negative speed, we therefore have the existence of a family of compactly supported sub-solutions for (\ref{e:vskew}).  Similarly, for any small positive speed we also have the existence of a super-solution for (\ref{e:vskew}).    Standard compactness arguments then give the existence of a stationary solution with limits $-1$ at $\xi=-\infty$ and $0$ at $\xi=+\infty$. The construction also gives $0<v_*<1$. We next show $v_*$ is monotonically decreasing. For this, note that for a minimum at $x_-$, $2u_*(x_-)-1+\mu\leq 0$ and for a maximum $x_+$, $2u_*(x_+)-1+\mu\geq 0$. Since $u_*$ is strictly decreasing,  $x_+\leq x_-$, a contradiction given that $v_*$ is decreasing near $\pm\infty$. 
 
Decay of this bounded solution in the $v$-equation is easily seen to be $\rme^{-(1+\sqrt{\mu})x}$, hence faster than the decay in the $u$-component, which establishes all conditions of Definition \ref{d:stab} provided we can show absence of unstable point spectrum in the linearization of the $v$-equation at the solution,
\[
v_t=v_{xx}+2v_x+(2u_*^2-1+\mu)(1-v_*^2)v = \mathcal{L}_v v.
\]
For this, we differentiate \eqref{e:vskew} at $v_*$ with respect to $x$ and find 
\[
\mathcal{L}_v v_{*,x}=h, \qquad h(x)=-4u_{*,x}u_* (v_*-v_*^3)>0.
\]
Suppose now that $\mathcal{L}_v \bar{v}=\lambda \bar{v}$ for some $\lambda>0$. Choosing the first eigenvalue, we may also assume that $\bar{v}(x)>0$ for all $x$. Inspecting decay at infinity, we see that $\bar{v}(x)$ decays faster than $-v_{*,x}$ so that there is a $\rho>0$ for which 
\begin{equation}\label{e:comp}
\rho\bar{v}(x)\leq -v_{*,x}(x) \text{, for all } x,\qquad \text{and } \qquad \rho\bar{v}(x_0)= -v_{*,x}(x_0) \text{, for some } x_0.
\end{equation}
At $x_0$ we then find, using \eqref{e:comp}, that $\rho\bar{v}(x_0)_x=-v_{*,xx}(x)$,  
\[
\mathcal{L}(-v_{*,x}-\rho\bar{v})(x_0)=-v_{*,xxx}(x_0)-\bar{v}_{xx}(x_0)=-h(x_0)-\lambda\rho \bar{v}(x_0).
\]
Using again \eqref{e:comp} we find $-v_{*,xxx}-\bar{v}_{xx}\geq 0$ but $-h-\lambda\bar{v}<0$, a contradiction. 

This shows the absence of unstable point spectrum and concludes the verification of all assumptions in Definition \ref{d:stab}. 

\paragraph{A conceptual existence proof.}
The fronts can also be obtained using a dynamical-systems type heteroclinic gluing argument. Such arguments were pursued in more difficult situations. In \cite{GS1}, an oscillatory invasion front follows a parameter step (rather than the primary front $u_*$, here), while in \cite{HSaccelerated} the secondary front is locked to the primary front due to a resonance pole in the linearization rather than the pinched double root.

\begin{figure}
\centering
\includegraphics[width=0.8\textwidth]{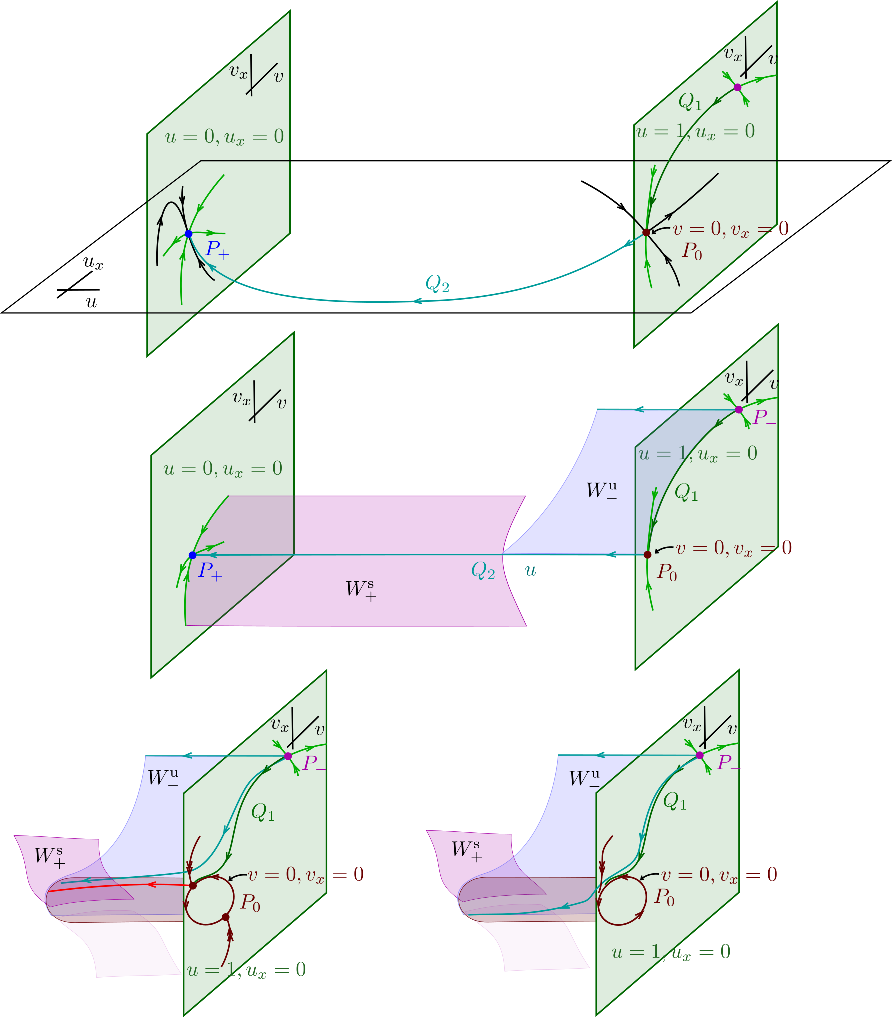}
\caption{\emph{Top:} Phase portrait in the four-dimensional space of $(u,u_x,v,v_x)$ here shown schematically, with the heteroclinics $Q_1$ and $Q_2$. \emph{Middle:} A simplified picture, replacing the fixed orbit $Q_2$ simply by a scalar heteroclinic on the line, showing now in more detail the dynamics in the $v$-$v_x$ directions.  One sees how $W_-^\mathrm{u}$ ``almost'' intersects $W^\mathrm{s}_+$ transversely. Upon increasing $\mu$, the double eigenvalue at $P_0$ becomes complex and $W_-^\mathrm{u}$ winds around the $v=v_x=0$ axis. \emph{Bottom:} This is made more visual in polar coordinates, where the $v$-$v_x$ axis is replaced by an invariant circle after choosing polar coordinates in this direction. The unstable manifold $W_-^\mathrm{u}$ is exponentially close to the center-unstable manifold $W^\mathrm{cu}_0$ of the saddle-node equilibrium on the invariant circle, which in turn intersects $W_+^\mathrm{s}$ transversely. For $\mu$ small, this intersection translates into an intersection of $W_-^\mathrm{u}$ and $W_+\mathrm{s}$ since trajectories in $W_-^\mathrm{u}$, fibered over the dynamics in $W^\mathrm{cu}_0$, are no longer blocked by the invariant strong unstable manifold of the saddle-node equilibrium.}
\label{f:het_gluing}
\end{figure}

Figure \ref{f:het_gluing} illustrates the heteroclinic bifurcation associated with the bifurcation at $\mu=0$. One glues two heteroclinic orbits: the first heteroclinic  $Q_1$ connects $P_-$ $u=v=1$, $u_x=v_x=0$ to $P_{0}$, $u=1,v=0$, $u_x=v_x=0$, within the invariant subspace $u=1,u_x=0$. The second heteroclinic connects $P_{0}$, $u=1,v=0$, $u_x=v_x=0$, to  $P_+$, $u=0,v=0$, $u_x=v_x=0$, within the invariant plane $v=v_x=0$. In order to find connections directly between $P_-$ and $P_+$, one tracks the unstable manifold of $P_-$ past a neighborhood of $P_0$ and matches with the stable manifold of $P_+$. Choosing polar coordinates in the $v$-$v_x$-plane, the double eigenvalue of the linearization at $P_0$ splitting into two complex eigenvalues as $\mu$ increases past $0$ in the $v$-$v_x$ plane corresponds to a saddle-node on the invariant circle representing the equilibrium. The unstable manifold of $P_-$ is exponentially close to the center-unstable manifold of the saddle-node equilibrium on the invariant circle and intersects the stable manifold of $P_+$ transversely. Increasing $\mu$ past $0$, this transverse intersection persists and now corresponds to an actual heteroclinic orbit. In fact, the picture also shows a family of heteroclinic orbits bifurcating at $\mu=0$, indexed by the number of turns around the circle $P_0$. The profiles exhibit an increasing number of sign changes and are expected to have an increasing number of unstable eigenvalues. We also refer to \cite{ADSS21} for background on the desingularization via polar coordinates and a similar analysis, used there to study fronts in bounded domains.

\begin{remark}[Distance scaling]\label{r:ds}
From the conceptual analysis, one finds immediately that the length of the plateau is determined by the passage time near a saddle-node bifurcation. In a scaled normal form, this saddle-node bifurcation is described as $\xi'=\mu+\xi^2$, and the time to move from $-\delta$ to $\delta$ in this equations scales with the inverse squareroot of $\mu$. The distance between fronts therefore scales as
\begin{equation}\label{e:distsqrt}
\Delta x\sim \frac{1}{\sqrt{\Delta\mu}}.
\end{equation}
This is in stark contrast to the front locking induced by a resonance pole, studied in detail in \cite{HSaccelerated}, where we found  much smaller separation distances,
\begin{equation}\label{e:distlog}
\Delta x\sim \log (\Delta\mu).
\end{equation}
\end{remark}

\paragraph{Conceptual stability.}
General results on heteroclinic gluing give that the linearization at the glued profile is the union of essential spectra at $P_-$ and $P_+$, eigenvalues at $Q_1$ and $Q_2$, and clusters of eigenvalues accumulating at the absolute spectrum of $P_0$ \cite{ssgluing}
. As a consequence, unstable eigenvalues can only bifurcate from the origin. Excluding these eigenvalues is somewhat subtle. In fact, the comparison arguments can be adapted to construct fronts with sign changes and establish the presence of unstable eigenvalues in the linearization near these non-monotone fronts. The fact that the first bifurcating front is in fact stable has however been shown in the more complicated case of an oscillatory secondary front in \cite{GdR}. We do not attempt to summarize the construction here.

\section{Analysis of Ginzburg-Landau with resonant forcing and Result \ref{res:3frontscgl}}

We start by considering real fronts and their stability. We then discuss secondary invasion front and evoke the picture from our skew-coupled toy example, Figure \ref{f:het_gluing}.
\begin{equation}\label{e:rcglf1}
\begin{aligned}
u_t&=u_{xx}+(1+\alpha) u +\beta (u^2-v^2)-u(u^2+v^2),\\
v_t&=v_{xx}+(1-\alpha) v +2\beta uv-v(u^2+v^2),
\end{aligned}
\end{equation}

\paragraph{Real fronts.}

In the real subsystem, we find the equation
\begin{equation}
u_t=u_{xx}+(1+\alpha)u +\beta u^2-u^3,
\end{equation}
with equilibria
\[
u=0,\qquad  u_\pm=\frac{1}{2}\left(\beta\pm \sqrt{4+4\alpha +\beta^2}\right).
\]
Fronts $u_*^\pm$ connecting $u_\pm$ to $u=0$ at the linear spreading speed $2\sqrt{1+\alpha}$ exist and are stable (that is, they are pulled) when $u_+/u_->-2$, which gives $\beta<\sqrt{(1+a)/2}$. This can be readily seen by finding explicit fronts where $u_x$ is a quadratic function of $u$ and using monotonicity in $c$.

\paragraph{Imaginary stability of real fronts.}
Linearizing the full equation \eqref{e:rcglf} at $u_*^\pm$, we find the linearized operator
\begin{equation}\label{e:rcglflin}
\mathcal{L}^\pm\begin{pmatrix}u\\v\end{pmatrix}=\begin{pmatrix}\mathcal{L}_u^\pm u\\\mathcal{L}_v^\pm v\end{pmatrix}=\begin{pmatrix}u_{xx}+c_\mathrm{lin}u_x+(1+\alpha) u +2\beta u_*^\pm u-3(u_*^\pm)^2u\\
v_{xx}+c_\mathrm{lin}v_x+(1-\alpha) v -2\beta u_*^\pm v-(u_*^\pm)^2 v
\end{pmatrix},
\end{equation}
which decouples, with $\mathcal{L}^\pm_u$ being simply the linearization at the real front, thus stable in an exponentially weighted norm. The linearization in the imaginary direction, $\mathcal{L}^\pm_v$ can exhibit instability.  We analyze the spectrum of $\mathcal{L}_v^\pm$ in three steps:

\emph{Leading edge.} In the leading edge, we can find the essential spectrum by setting $u_*=0$. In an exponential weight $v(x)\rme^{\eta x}$, $\eta=\sqrt{1+\alpha}$, determined by the real operator, the essential spectrum can be found after Fourier transform with Fourier exponent $k=\tilde{k}+\rmi\eta$, $\tilde{k}\in\R$. A quick calculation finds that the essential spectrum is stable for $\alpha>0$  with a maximum real part at $-2\alpha$, for both  $\mathcal{L}_v^\pm$.

\emph{Wake.} In the wake we set $u_*=u_*^\pm$ and find the associated dispersion relation 
\[
d^v_\pm(\lambda,\nu)=\nu^2+2\sqrt{1+\alpha}\nu +1-\alpha -2\beta u_*^\pm -(u_*^\pm)^2 -\lambda=0,
\]
Without exponential weights, stability is determined by setting $\nu=0$. One readily finds that, at $u_*^+$, $\lambda\leq -2\alpha-\frac{3}{2}\beta(\beta+\sqrt{4+4\alpha +\beta^2})<0$, hence stability. At $u_*^-$, we find a maximal real part as 
\[
\lambda= -2\alpha+\frac{3}{2}\beta(-\beta+\sqrt{4+4\alpha +\beta^2}),
\]
which is unstable for $\beta>\frac{2\alpha}{\sqrt{9+3\alpha}}$. Similarly, in the optimal exponential weight where with $\nu=-\sqrt{1+\alpha}$, we find instability when $\beta>\beta_\mathrm{lin}=\frac{1+3\alpha}{\sqrt{6}}$.

\begin{figure}
\centering
\includegraphics[width=0.43\textwidth]{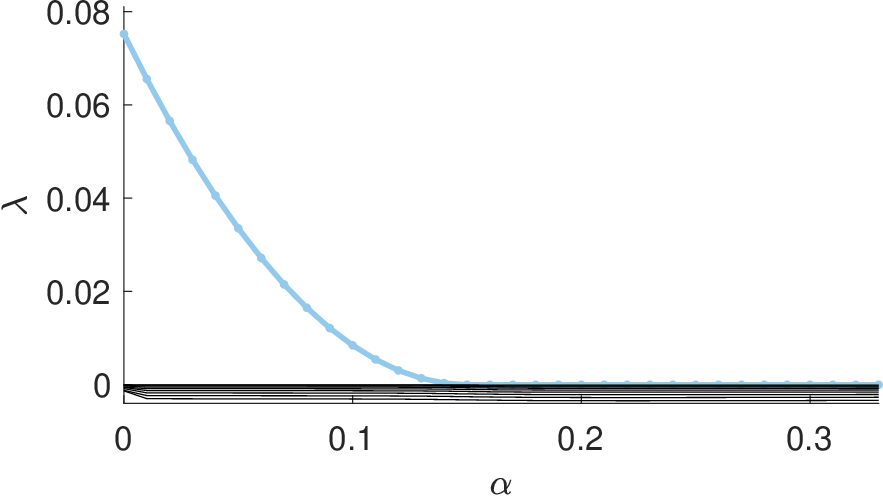}\qquad
\includegraphics[width=0.43\textwidth]{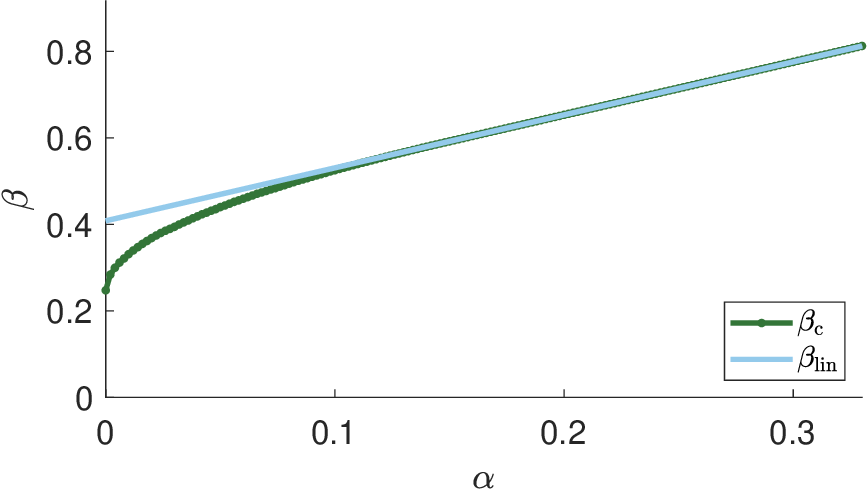}

\caption{\emph{Left:} Spectra of $\tilde{\mathcal{L}}_v^-$ at  $\beta=\frac{1+3\alpha}{\sqrt{6}}$, computed numerically with numerically computed profile $u_*$ as the critical front at speed $2\sqrt{1+\alpha}$ in the $u$-equation. For $\alpha>\alpha_\mathrm{c}\sim 0.15$, there are no unstable eigenvalues so that increasing $\beta$ past this threshold produces an instability caused by the essential spectrum only. For $\alpha<\alpha_\mathrm{c}$, the instability occurs at a value $\beta=\beta_\mathrm{c}<\frac{1+3\alpha}{\sqrt{6}}$ due to point spectrum and leads to an alternative form of locking between primary and secondary invasion, discussed in \cite{HSaccelerated}. \emph{Right:} Critical value of $\beta$ for which the linearization $\tilde{\mathcal{L}}_v^-$ is neutrally stable. For values above this critical $\beta$, there exist locked fronts. The blue curve shows the prediction from the essential spectrum in the wake, the green curve shows the earlier locking due to an eigenvalue. }\label{f:schsp}
\end{figure}

\emph{Intermediate regime.} 
We now consider eigenvalues in the optimal exponential weight $\eta=-\sqrt{1+\alpha}$, which leads to studying the operator
\[
\tilde{\mathcal{L}}_v^\pm v=v_{xx}-2\alpha v -2\beta u_*^\pm v-(u_*^\pm)^2 v.
\]
Clearly, $\tilde{\mathcal{L}}_v^+$ is a self-adjoint Schr\"{o}dinger operator and negative semi-definite since $u_*^+>0$. It turns out that $\tilde{\mathcal{L}}_-$  may have positive eigenvalues. We computed the spectrum
at the critical value $\beta=\beta_\mathrm{lin}=\frac{1+3\alpha}{\sqrt{6}}$, in which case the essential spectrum touches the origin; see Fig. \ref{f:schsp}. For $\alpha<\alpha_\mathrm{c}\sim 0.15$, $
\tilde{\mathcal{L}}_v^-$ possesses a positive eigenvalue at $\beta=\beta_\mathrm{lin}$. It is negative semi-definite for $\beta\leq \beta_\mathrm{c}<\beta_\mathrm{lin}$. Setting $\beta_\mathrm{c}=\beta_\mathrm{lin}$ for $\alpha>\alpha_\mathrm{c}$ then gives instability in the regime described in Result \ref{res:3frontscgl}.

\paragraph{Secondary fronts.}

It appears difficult to establish existence and stability of pulled secondary fronts, that is, fronts invading $A_3$ and leaving behind $A_2$ (or $A_4$ for the complex conjugate front). Nevertheless, these fronts are easy to observe in direct simulations and can be computed numerically, both at the predicted linear invasion speed. We investigated these fronts numerically, in direct simulations. Sweeping the parameter plane with $0<\alpha<1/3$ and $\beta_\mathrm{p}/2<\beta<\beta_\mathrm{p}$, we found invasion of $A_3$ by $A_2$ at the linear spreading speed. In particular, we did not find pushed fronts in this invasion process.

\paragraph{Gluing.}

\begin{figure}
\begin{minipage}{0.33\textwidth}\includegraphics[width=\textwidth]{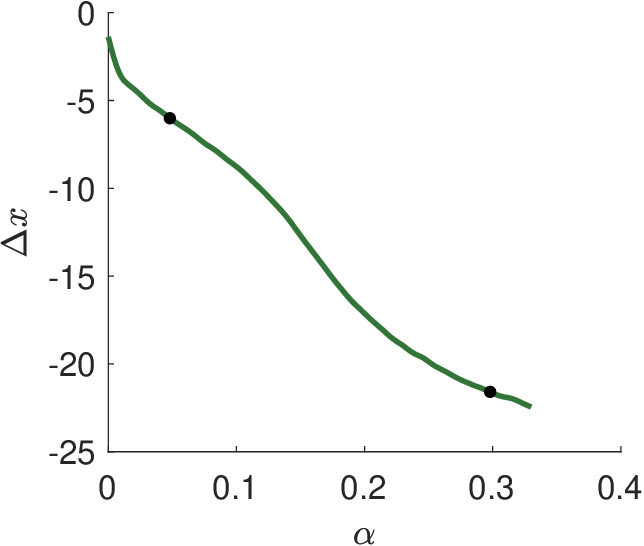}
\end{minipage}\hfill
\begin{minipage}{0.25\textwidth}
\includegraphics[height=1in]{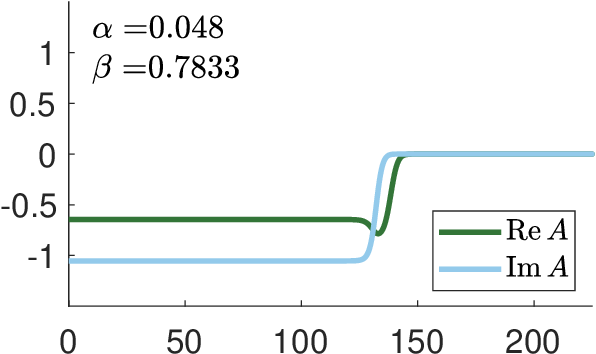}\\[0.1in]
\includegraphics[height=1in]{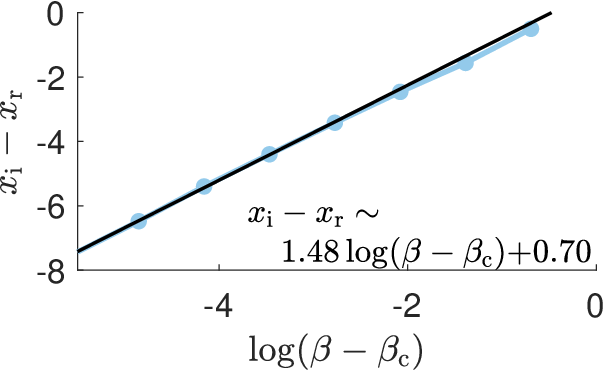}
\end{minipage}\hfill
\begin{minipage}{0.25\textwidth}
\includegraphics[height=1in]{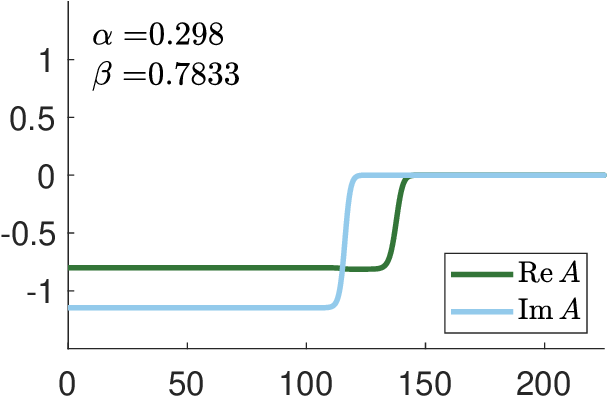}\\[0.1in]
\includegraphics[height=1in]{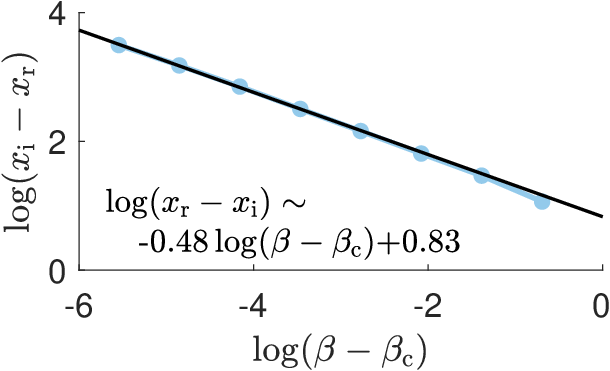}
\end{minipage}
\caption{\emph{Left:} Location of secondary front relative to primary front, measured by the difference of locations $x_\mathrm{i}-x_\mathrm{r}$ where $\Im A(x_\mathrm{i})=0.1$ and $\Re A(x_\mathrm{r})=0.1$. The distance increases in $\alpha$.
One also notices a qualitative change between front splitting induced by a resonance pole, $\alpha<\alpha_\mathrm{c}$, and front splitting induced by essential spectrum,  $\alpha>\alpha_\mathrm{c}$, $\alpha_\mathrm{c}\sim0.15$; see ; see \eqref{e:distsqrt} and \eqref{e:distlog}. The former leads to distances growing with $\log(\beta-\beta_\mathrm{c}$ and the latter to distances growing with $(\beta-\beta_\mathrm{c})^{-1/2}$.
\emph{Center and Right:} Sample profiles (top row) and distance as a function of $\beta$ for fixed $\alpha$ (bottom row), demonstrating the $\log$-scaling below $\alpha_\mathrm{c}$ and the inverse square root scaling above $\alpha_\mathrm{c}$.
} \label{f:glue}
\end{figure}

We also studied the glued profiles numerically; see Figure \ref{f:glue}. We found fronts where $A_2$ and, by symmetry, $A_4$ invade the origin for all $\beta\in(\beta_\mathrm{c},\beta_\mathrm{p})$, the range described in Result \ref{res:3frontscgl}. We also found that the fronts indeed disappear upon decreasing $\beta$ past $\beta_\mathrm{c}$ by splitting into two fronts, with first $A_3$ invading the origin, followed by a slower invasion of $A_3$ by $A_{2/4}$. The distance between the two fronts diverges as $\beta\searrow \beta_\mathrm{c}$. Figure \ref{f:glue} shows the distance and front profiles at $\beta_\mathrm{c}+0.01$.

We also note that one expects a long-range interaction between primary and secondary front for $\beta<\beta_\mathrm{c}$, $\alpha<\alpha<\mathrm{c}$: the resonance pole in the linearization of the primary front leads to an accelerated speed of the secondary front although the distance increases linearly in time; see \cite{HSaccelerated}.

\section{Discussion}\label{s:dis}

We established rigorously that reaction-diffusion equations can possess more than two selected pulled fronts. This contradicts the simple intuition that pulled fronts are selected by their Gaussian tail in the leading edge, which is either positive or negative in the direction of the eigenvector. We construct examples in scalar reaction-diffusion equations and in a system of skew-coupled equations. In both cases, the fronts arise through the concatenation of a primary and a secondary front. In the first example, the intermediate state is stable, in the latter it is unstable. We also study, partly numerically, an example that arises as a modulation equation near a Turing instability in a spatially periodic medium, where multiple fronts again arise through concatenation, here at an unstable state. Our analysis emphasizes the robustness of the phenomena seen here, showing in particular that the higher multiplicity is not due to a degeneracy, as for instance the skew-coupling, or a degenerate leading edge, as in the complex Ginzburg-Landau equation. We therefore rely on the conceptual approach towards front selection based on marginal stability initiated for linear equations in \cite{HolzerScheelPointwiseGrowth} and for nonlinear equations in \cite{as1,avery2}: we show existence and stability and conclude selection and robustness using the conceptual results in \cite{as1,avery2}.

In our example, we find pulled fronts but we suspect that the analysis can readily be adapted to find pushed fronts, with similar multiplicities. In fact, when the strength $\beta$ of quadratic terms is larger, one expects pushed fronts.

In the pattern-formation example, it would be interesting to derive a more comprehensive picture of the existence and bifurcation of pulled and pushed fronts in the parameters $\alpha,\beta$, possibly describing in more detail the boundaries of the basin of attraction.  We emphasize however that such basin boundaries are poorly understood even in the example of the Nagumo equation \eqref{e:nag}.

Returning to the terminology of a pulled front, one may wonder how then exactly the selection of a state in the wake occurs. In an exponentially weighted space, the leading edge dynamics are simply diffusive, induced by the expansion of the dispersion relation near the double root $\lambda\sim d_\mathrm{eff}\nu^2$.
From this perspective, the leading order dynamics \emph{do not} exhibit directional transport! This is different for fronts with speeds $c>c_\mathrm{lin}$, for which the dispersion relation starts with $\lambda\sim (c-c_\mathrm{lin})\nu+d_\mathrm{eff}\nu^2$ 
in an appropriately weighted space, and which therefore possess directed transport with group velocity $c_\mathrm{g}=c-c_\mathrm{lin}$ 
from the leading edge to the wake. 

\begin{figure}
\centering
\includegraphics[width=0.33\textwidth]{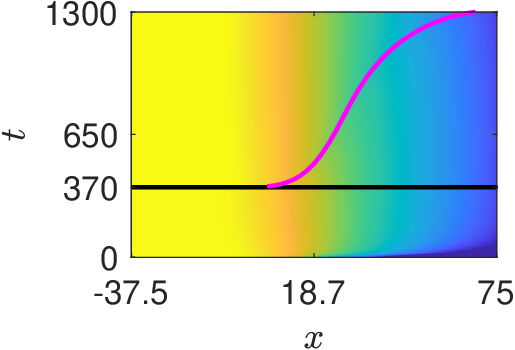}\hspace*{0.7in}
\includegraphics[width=0.33\textwidth]{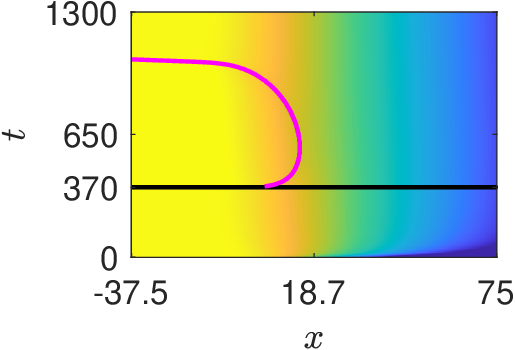}
\caption{Simulation of \eqref{e:cglf} with $\alpha=0.02$ and $\beta=0.7$, initial conditions $A_0(x)=\rme^{4\pi\rmi/3}$ for $x<0$, zero elsewhere, in $x\in[-150,75]$ (only part of the domain is shown). After time $t=370$, the solution is replaced by $A=1$ in $x<L_0$, unchanged in $x>L_0$, with $L_0=-1.82$ (left) and $L_0=-2.3$ (right). Shown is a plot of $\log(|\Re A|)$, the time when the solution is altered, and in magenta a space-time curve where $\Re A$ changes sign. The sign change from the wake propagates to the leading edge when a large enough part of the solution is altered (left) and recedes to recover the original state in the wake when the change is too far in the wake (right). }
\label{f:mismatch}
\end{figure}

For the critical front at speed $c_\mathrm{lin}$, selection is established from steep initial conditions by matching a Gaussian profile in the leading edge with the nonlinear front in the wake. Errors in the matching problem,  such as a mismatch from gluing tail and wake from different fronts,  are controlled using the linearization, and may, in the absence of transport lead to corrections in either the leading edge or the tail. We did in fact observe both possibilities as shown in Figure \ref{f:mismatch}.

In a different direction, one could wonder about oscillatory fronts, when linear marginal stability is caused by oscillations in the leading edge, $d_c(\rmi\omega_*,\nu_*)=0$, $\partial_\nu(\rmi\omega_*,\nu_*)=0$. These situations arise in particular near pattern-forming instabilities or in phase-separation models \cite{CE1990,S2017}.  In fact, the Ginzburg-Landau equation \eqref{e:cgl} is the modulation equation near such an instability and the double double root at the origin corresponds to two double roots at $\pm\rmi\omega_*$ in a pattern-forming equation such as a reaction-diffusion system near a Turing instability. Other examples include in particular the Ginzburg-Landau equation with complex coefficients. From our perspective, the linear equation possesses a unique complex eigenvector, which one might expect to be reflected in a uniqueness of selected invasion fronts,  of course up to time shifts, which correspond to complex phase rotation of the eigenvector, and to complex phase rotation of fronts in Ginzburg-Landau models. We suspect that the constructions here can be adapted to give examples of non-uniqueness of selected oscillatory invasion fronts as well.
%
%

%
\end{document}